\documentclass{aims} 
\usepackage{algorithm}
\usepackage{comment}
\usepackage{amsmath}
\usepackage{algorithmic}
\usepackage{multirow}
\usepackage{longtable}
\usepackage{array}
  \usepackage{paralist}
  \usepackage{graphics} 
  \usepackage{epsfig} 
\usepackage{graphicx}  
\usepackage[colorlinks=true]{hyperref}
\hypersetup{urlcolor=blue, citecolor=red}

\usepackage{bm}
\allowdisplaybreaks

\def\currentvolume{X}
 \def\currentissue{X}
  \def\currentyear{20XX}
   \def\currentmonth{XX}
    \def\ppages{X--XX}
     \def\DOI{10.3934/xx.xxx}

\newtheorem{theorem}{Theorem}[section]

\newtheorem{proposition}{Proposition}

\theoremstyle{definition}
\newtheorem{definition}[theorem]{Definition}
\newtheorem{remark}{Remark}

\title[permutation XOR cellular automata]
{Permutation XOR cellular automata: partial-shift connections and half-primitive polynomials
}

\author[Suzuki, Okuyama, Takahashi, Saito]{}


\thanks{$^*$  Corresponding author: Toshimichi Saito}

\begin{document}
\maketitle

\centerline{\scshape Yosuke Suzuki,  Naoto Okuyama, Hiroki Takahashi and Toshimichi Saito}
\medskip
{\footnotesize
 \centerline{Department of Electrical and Electronic Engineering}
   \centerline{HOSEI University, Tokyo, 184-8584, Japan}
} 

%

\bigskip



\begin{abstract}
This paper presents a permutation XOR cellular automaton, a simple digital dynamical system composed of 
an elementary cellular automaton of the XOR rule and a partial-shift permutation connection characterized by a shift-part parameter.  
The proposed system is suitable for precise analysis and simple hardware implementation. 
Depending on the partial-shift connection and system dimension, the automaton can generate long binary periodic orbits (LBPOs) with  strong stability. 
As a principal result, we identify the partial-shift parameters and dimensions that enable the generation of LBPOs. 
In our analysis, these LBPOs are characterize using half-primitive polynomials over $\mathrm{GF}(2)$. 
Furthermore, we present a simple FPGA-based hardware implementation. 
The hardware transforms typical LBPOs into electrical signals for potential engineering applications. 
\end{abstract}

\section{Introduction}
\label{intro}
An elementary cellular automaton (ECA) is a simple digital system with rich dynamical behavior \cite{ca1} \cite{ca2} \cite{ca3} \cite{ca4} \cite{ca5} \cite{ca6}. 
Its dynamics is governed by a local update rule defined by a simple Boolean function and is described by an autonomous difference equation for binary state vectors. 
Depending on the rule, an ECA can generate a wide variety of spatiotemporal patterns. 
Since the number of possible binary states is finite, every sequence eventually falls into a periodic orbit of binary vectors (BPO). 
Real/potential applications of the BPOs include
reservoir computing \cite{rc1} \cite{rc2} \cite{rc3} \cite{rc4}, 
error-correcting codes \cite{error}, and 
data description \cite{data}. 
Therefore, understanding the dynamics of ECAs is important from both fundamental and applied perspectives. 

Inspired by ECAs, we have introduced the permutation elementary cellular automaton (PECA, \cite{iconip22} \cite{iconip24} \cite{ncom25})
which is obtained by incorporating a one-to-one permutation connection into an ECA. 
Through exhaustive numerical analysis of low-dimensional PECAs, 
we have shown that the permutation connection enables the generation of various BPOs that cannot occur in ECAs \cite{iconip22} \cite{taka} \cite{mikito}. 
However, the number of possible permutation connections grows 
factorially with the system dimension, making the analysis computationally intractable because of the curse of dimensionality. 

To analyze higher-dimensional BPOs, this paper presents a permutation XOR cellular automaton (PXCA), which can be regarded as a simplified version of the PECA characterized by two key features. 
First, the permutation connections are restricted to a partial-shift connections characterized by the shift-part parameter $r$. 
This constitutes the main novel idea of this paper. 
As the dimension $N$ increases, the number of possible partial-shift connections grows only linearly, thereby we can escape from the curse of dimensionality. 
Second, the ECA rules are restricted to a simple XOR rule. 
Both the partial-shift connection and the XOR rule facilitate precise analysis of the dynamics as well as efficient hardware implementation.  

Depending on the parameter $r$ and the dimension $N$, the PXCA can generate a wide variety of BPOs. 
Since complete analysis of all BPOs is difficult, we focus on long binary periodic orbits (LBPOs) with strong stability. 
Such stability is related to both size of the basin of attraction \cite{ca3} \cite{ca4} and error-correcting capability. 
To analyze the LBPOs, we introduce the transition matrix of the PXCA and derive its characteristic polynomial $q(x)$ over the Galois field $\mathrm{GF}(2)$ \cite{lfsr1} \cite{lfsr2}. 
We show that an LBPO is characterized by an $N$-th degree half-primitive polynomial $q(x) = xp(x)$, where $p(x)$ is primitive. 
As a main result, through precise analysis, we identify the values of $r$ and $N$ that enable the generation of LBPOs. 
 
As a fundamental step to engineering applications of PXCAs, we present a simple hardware implementation on a Field-Programmable Gate Array (FPGA). 
The FPGA-based hardware converts LBPOs into electrical signals. 
The electrical signals can be used in application systems such as pseudo-random number generators \cite{rnd} and control of switching circuits \cite{pe1} \cite{pe2}. 
Owing to the strong stability of LBPOs, 
the generated signals exhibit error-correction capability and robustness against external perturbations. 

As a remark on the novelty of this work, our previous study \cite{ncom25} investigated a PXCA with a full-shift permutation structure. 
However, it did not address partial-shift connections, LBPOs, or the characteristic polynomial analysis presented in this paper. 

\section{Preliminary}

\begin{figure}[b]
\centering
\includegraphics[width=0.5\columnwidth]{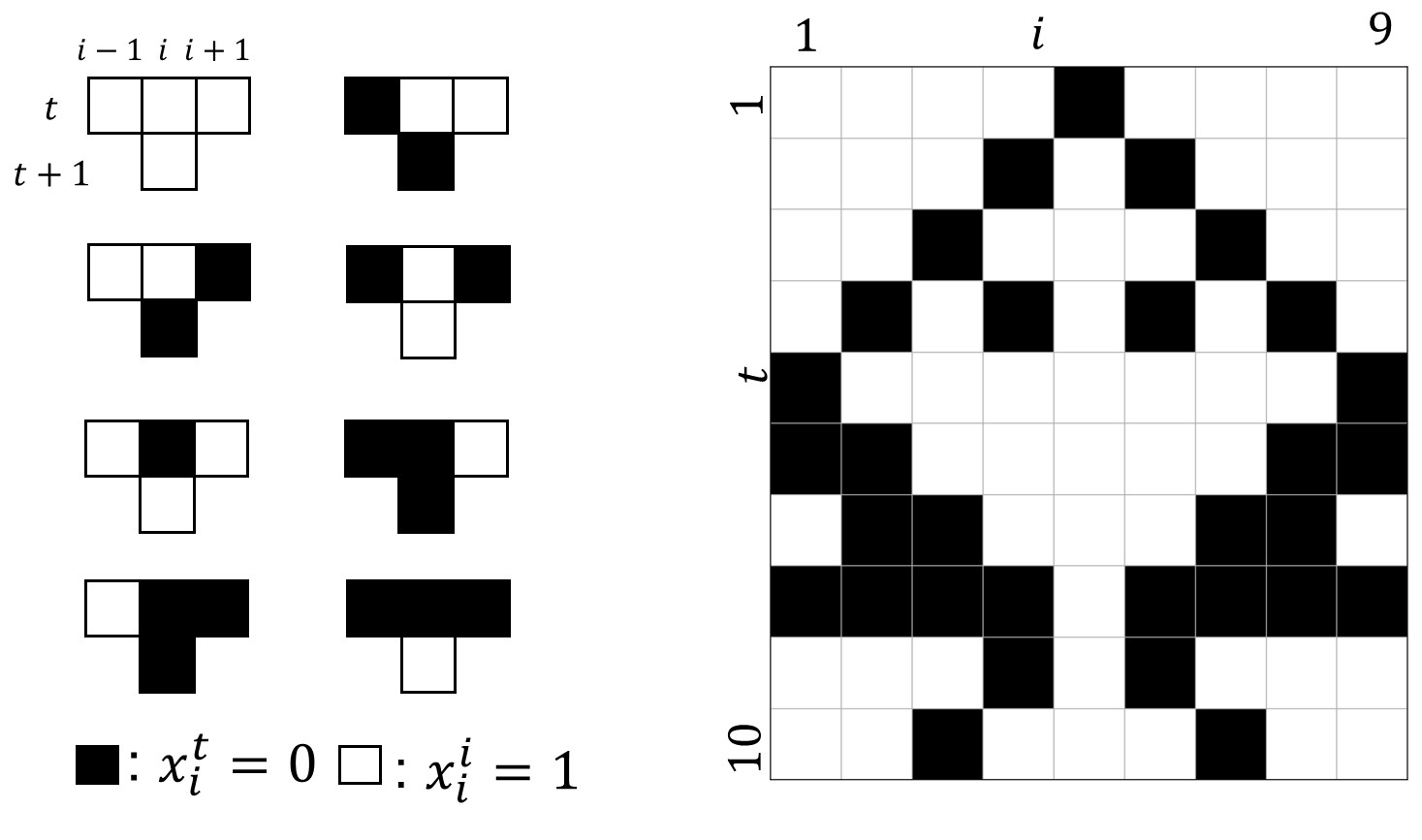}
\caption{ECA rule 90 and a spatiotemporal pattern}
\label{fg1}

\vspace*{5mm}

\centering
\includegraphics[width=0.35\columnwidth]{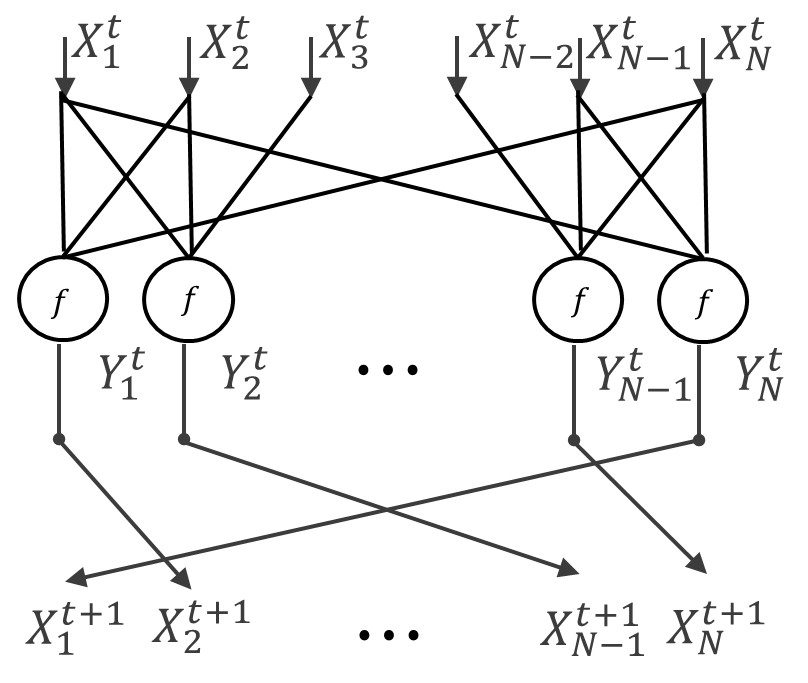}
\caption{PECA. $f$ denotes a Boolean function of 3 inputs. }
\label{fg2}
\end{figure}

As a preliminary step, we introduce two  fundamental dynamical systems: the elementary cellular automaton (ECA \cite{ca1}) and the permutation elementary cellular automaton (PECA \cite{iconip22}). 
An ECA is defined by the following autonomous difference equation for binary state variables: 
\begin{equation}
X_i^{t+1}=f(X^t_{i-1}, X^t_{i}, X^t_{i+1}), \ i \in \{1, \cdots, N \}, \ N \ge 5, 
\end{equation}
where $X^t_i \in \{0, 1\}$ denotes the $i$-th binary state variable at discrete time $t$. 
The state variables $X_i$ are arranged on a ring, i.e., $X^t_{N+1} \equiv X^t_1$ and $X^t_0 \equiv X^t_N$. 
The dynamics is governed by a rule of Boolean function $f$ from three binary inputs to one binary output. For example, 
\begin{equation}
\begin{array}{cccc}
f(1,1,1)=0 \ & \ f(1,1,0)=1 \ & \ f(1,0,1)=0 \ &  \ f(1,0,0)=1\\
f(0,1,1)=1 \ & \ f(0,1,0)=0 \ & \ f(0,0,1)=1 \ & \ f(0,0,0)=0
\end{array}
\label{rule90}
\end{equation}
A rule is identified by the decimal representation of its 8 output bits. 
Since 
$(01011010)_2=(90)_{10}$, the rule defined in Eq. \eqref{rule90} is called Rule 90. 
The total number of possible rules is $2^8=256$. 
Depending on the rule, an ECA can generate a wide variety of binary-state sequences. 
For example, Rule 90 generates the sequence shown in Fig.\ref{fg1}. 
Since the number of possible binary states is finite, every sequence eventually falls into a periodic orbit of binary vectors (BPO). 

Incorporating a one-to-one permutation connection into an ECA yields a PECA. 
As illustrated by the three-layer architecture in Fig.\ref{fg2}, the input-to-hidden layer implements an ECA whereas the hidden-to-output layer implements a permutation connection. 
The dynamics is described by
\begin{equation}
    \begin{array}{l}
        X^{t+1}_i = Y_{\sigma(i)}^{t},\
Y_i^{t}=f(X^t_{i-1}, X^t_{i}, X^t_{i+1}), \ \ 
        \sigma =
        \left(
        \begin{array}{cccc}
            1 & 2 & \cdots & N\\
             \sigma(1)& \sigma(2) & \cdots & \sigma(N)
        \end{array}
        \right) \\

    \end{array}
\label{peca}
\end{equation}
where $X^t_i \in \{0, 1\}$ denotes the $i$-th binary state variable,  $Y_i^t$ denotes the $i$-th binary hidden state variable, and $\sigma$ is a permutation operator. 
Through exhaustive numerical analysis of low-dimensional PECAs, we have shown that the permutation connection $\sigma$ enables the generation of various BPOs that cannot occur in the ECAs \cite{iconip22}. 
However, as the system dimension $N$ increases, the number of possible permutation connections grows factorially, making exhaustive analysis computationally infeasible because of the curse of dimensionality. 

\section{Permutation XOR Cellular Automaton} 

\subsection{Definition and examples}

\begin{figure}[b]
\centering
\includegraphics[width=0.8\columnwidth]{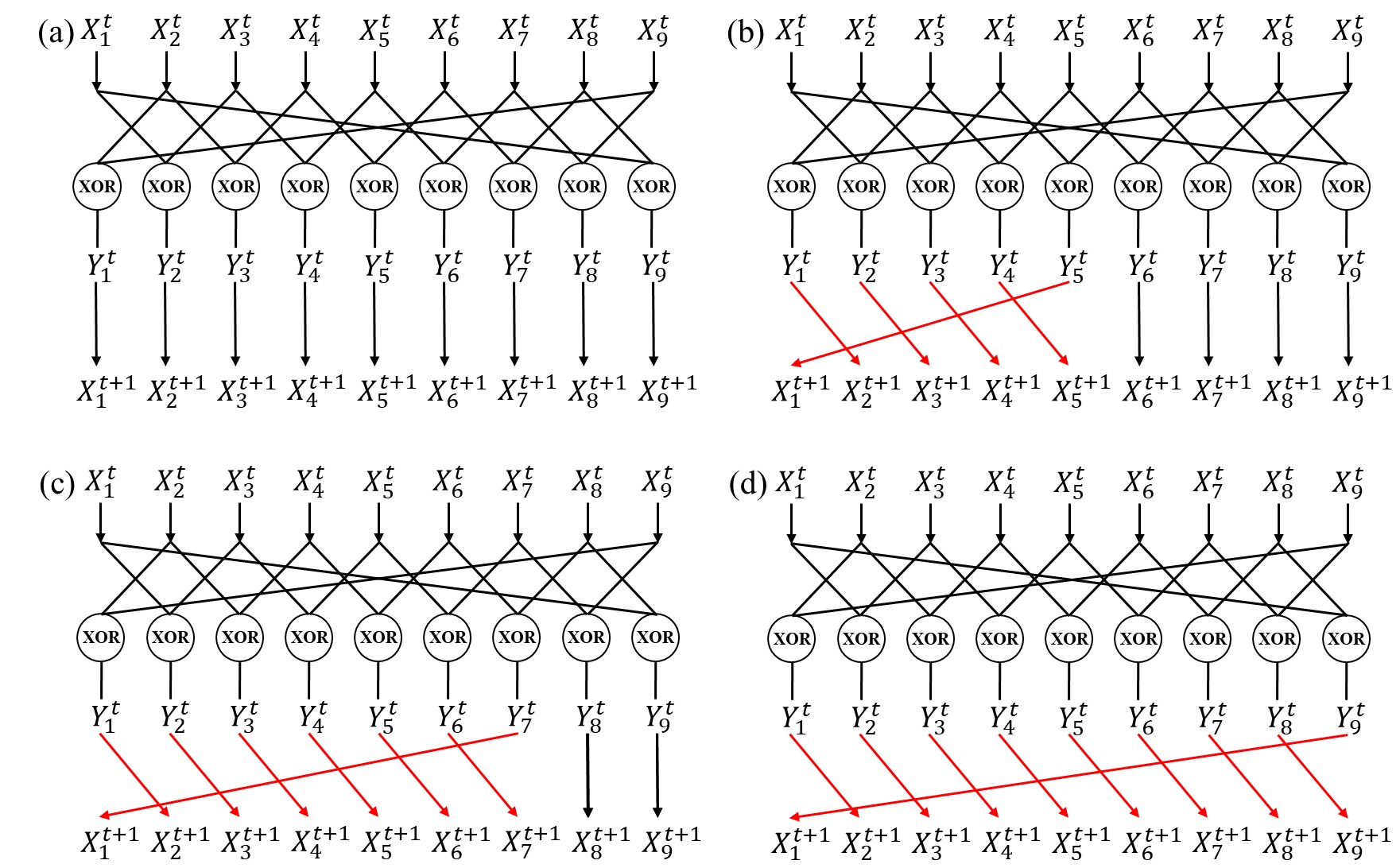}
\caption{PXCA with partial-shift connection (red). (a) $r=1$ (ECA). (b) $r=5$. (c) $r=7$. (d) $r=9$ (full shift). }
\label{fg3}
\end{figure}

\begin{figure}[tb]
\centering
\includegraphics[width=0.7\columnwidth,]{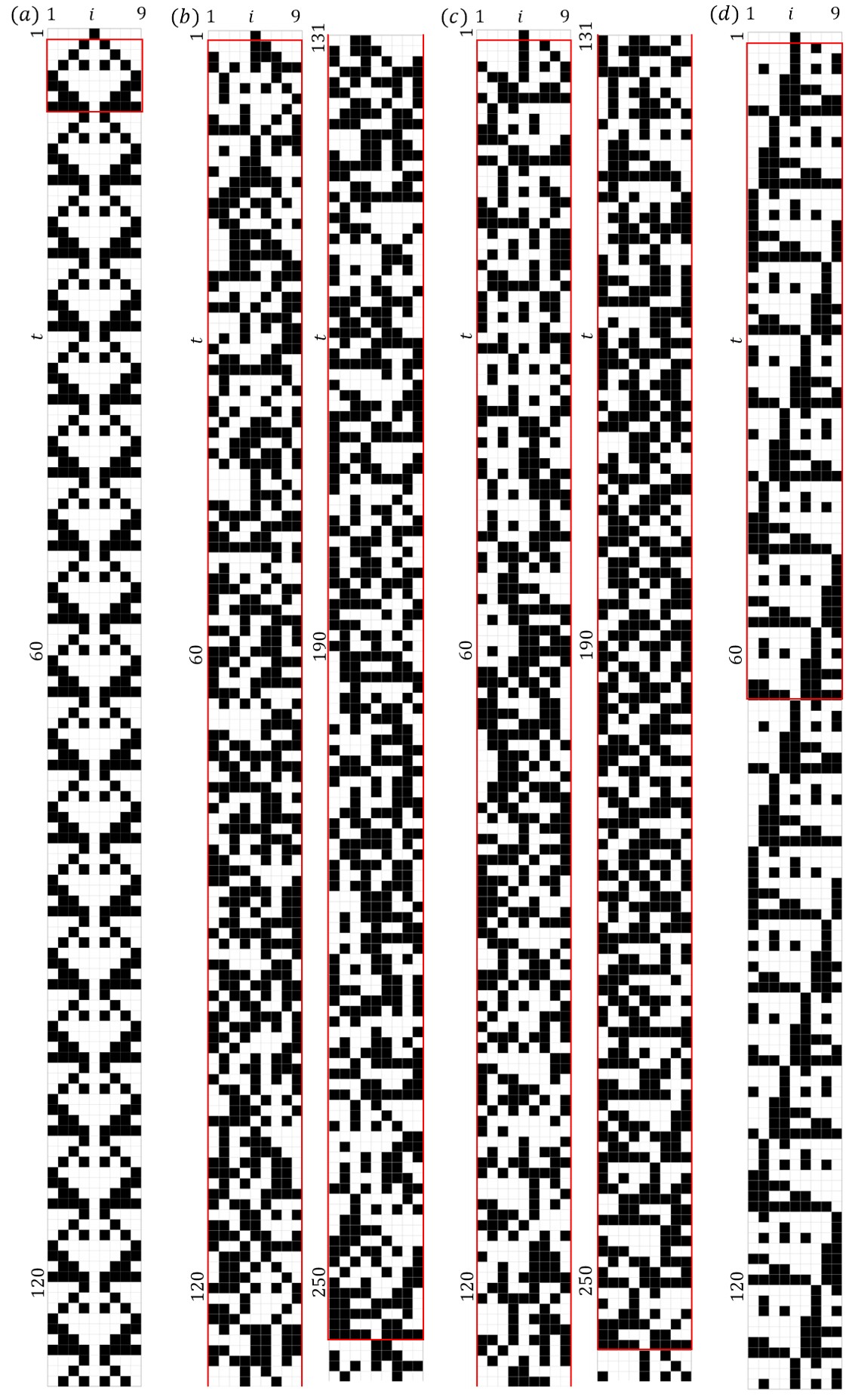}
\caption{BPOs from PXCA in Fig. \ref{fg3}. (a) Period 7 for $r=1$.  (b) and (c) Period 255 for $r=5$ and $r=7$. (d) Period 63 for $r=9$. }
\label{fg4}
\end{figure}

In order to enable efficient analysis, 
we present the permutation XOR cellular automaton (PXCA), a simplified version of the PECA with two key features. 
{\it First, the permutation connections are restricted to a partial-shift connection, denoted by} 
\[
\sigma_r = 
\left(\begin{array}{ccccccc}
1 & 2 & \cdots &  r  &  r+1 & \cdots & N\\
r & 1 & \cdots & r-1 &  r+1 & \cdots & N
\end{array}\right) \equiv 
        \left(
        \begin{array}{cccc}
            1 & 2 & \cdots & N\\
             \sigma_r(1)& \sigma_r(2) & \cdots & \sigma_r(N)
        \end{array}
        \right)
\]
The partial-shift connection is characterized by the shift-part parameter $r \in \{2, \cdots, N \}$. 
For $r=1$, let $\sigma_r$ be the identity: the PXCA becomes equivalent to an ECA. 
This restriction constitutes the main novel idea of this paper. 
Fig. \ref{fg3} shows several examples for $N=9$. 
Since the number of possible partial-shift permutations grows only linearly with $N$, we can escape from the curse of dimensionality.  
Second, we restrict the ECA rule to Rule 90 (Eq.~\eqref{rule90}), defined by the XOR operation:  
\[
f(X_{i-1}, X_{i}, X_{i+1}) = X_{i-1} \oplus X_{i+1}
\]
Rule 90 (equivalent to Rule 165 under state inversion) is known to exhibit a variety of interesting dynamical behaviors \cite{ca1}.  
Thus, the PXCA dynamics is described by
\begin{equation}
    \begin{array}{l}
        X^{t+1}_i = Y_{\sigma_r(i)}^{t+1},\
        Y_i^{t+1}=X^t_{i-1} \oplus X^t_{i+1}, \ \ 
        \sigma_r(i) =\left\{
        \begin{array}{ll}
         r & \text{for } i = 1 \\
        i-1 & \text{for } 2 \le i \le r\\
            i & \text{for } r+1 \le i \le N
        \end{array}\right.
    \end{array}
\label{pxca}
\end{equation}
%
The PXCA, consisting of the partial-shift connection and Rule 90, facilitates both precise analysis of the dynamics and efficient hardware implementation. 
For convenience, we introduce the vector form of Eq. \eqref{pxca}
\begin{equation}
    \begin{array}{l}
    \bm{X}^{t+1}=F(\bm{X}^t), \ \  \bm{X}^t \equiv (X^t_1, \cdots, X^t_{N}) \in \bm{B}^N
    \end{array}
    \label{vform}
\end{equation}
where $\bm{B}^N$ denotes the set of all $N$-dimensional binary vectors. 
Depending on the shift-part parameter $r$ and the dimension $N$, the PXCA exhibits a wide variety of BPOs. 
Fig.~\ref{fg3} shows several examples of 9-dimensional PXCAs ($N=9$), where the input-to-hidden layer is governed by Rule 90 and the hidden-to-output layer is implemented by the partial-shift connection. 
We can see the following. 
\begin{itemize}
\item (a) $r=1$ (Rule 90): 
the PXCA generates a BPO with period 7 as shown in Fig.~\ref{fg4}(a). 
This PXCA has multiple BPOs \cite{ncom25} and generates one of them depending on the initial condition. 
\item (b) $r=5$: the PXCA generates the BPO with period 255 in Fig.~\ref{fg4}(b). 
This is a long binary periodic orbit (LBPO)  discussed in the next section. 
\item (c) $r=7$: the PXCA generates another LBPO in Fig.~\ref{fg4}(c): the period is the same, but the sequence is different. 
The PXCA generates a LBPO iff $r \in \{5, 7\}$.
\item (d) $r=9$ (full-shift permutation): 
the PXCA generates a BPO with period 63 in Fig.~\ref{fg4} (d). 
This PXCA has multiple BPOs \cite{ncom25}.  
\end{itemize}

\subsection{Long Binary Periodic Orbit}
First, we give a fundamental definition. 

{\definition 
A binary vector $\bm{X}_p$ is called a binary periodic point (BPP) with period $p$ if it returns to itself after $p$ iterations: $\bm{X}_p = F^p(\bm{X}_p)$, 
$\bm{X}_p\neq F^k(\bm{X}_p)$ for $ 1 \le k < p$, where $F^k$ is the $k$-fold composition of $F$. 
A sequence of BPPs, $( F(\bm{X}_p), F^2(\bm{X}_p), \dots, F^{p}(\bm{X}_p) )$,  
is called a binary periodic orbit (BPO). 
A binary vector $\bm{X}_d$ is called a direct eventually periodic point (DEPP) of a BPO if it is not a BPP but falls directly into the BPO: $\bm{X}_d \ne$ BPP and 
$F(\bm{X}_d) \in$BPO.  
}

\begin{figure}[tb]
\centering
\includegraphics[width=0.7
\columnwidth]{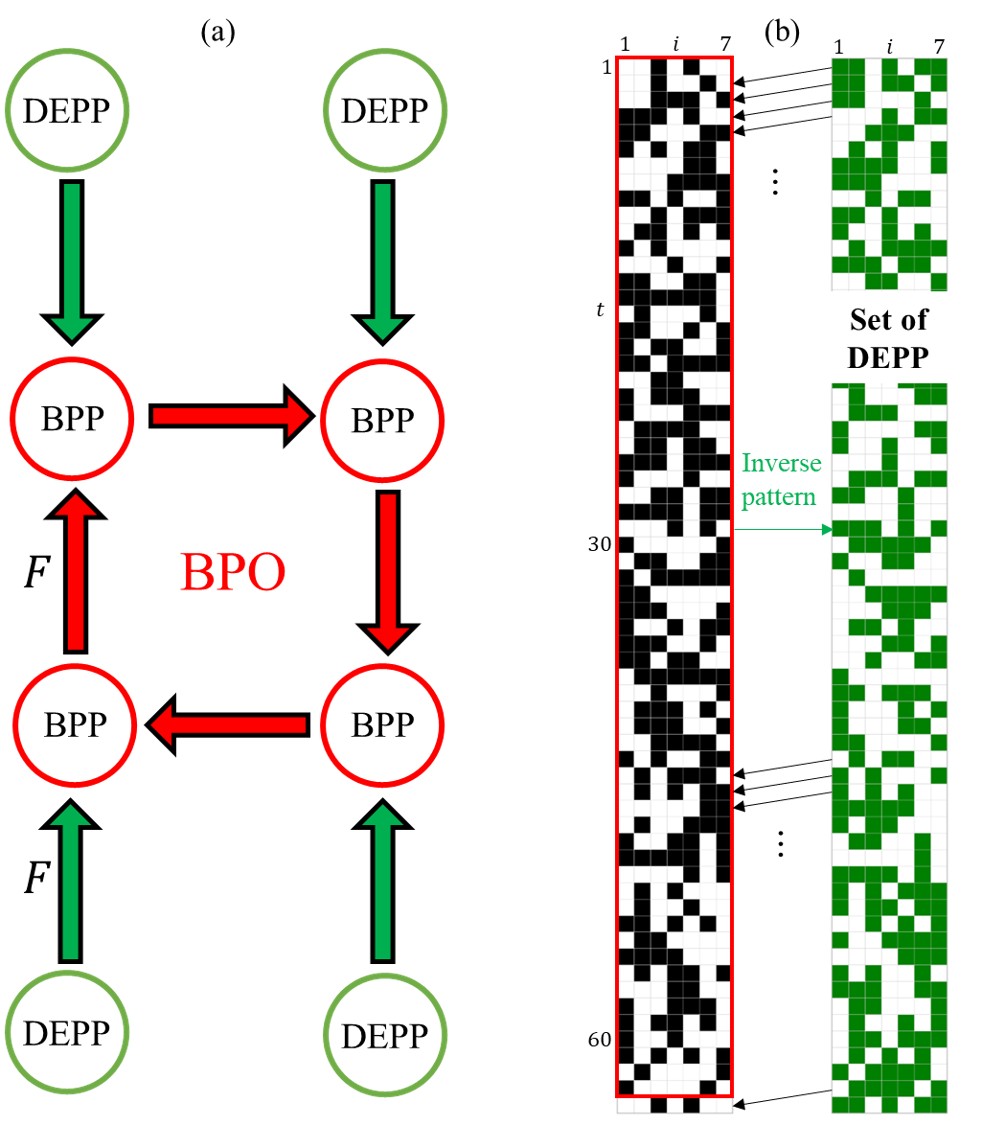}
\caption{BPO, BPP, and DEPP. 
(a) Conceptual diagram.\\ 
(b) LBPO for $N=7$ and $r=3$. 
}
\label{fg5}
\end{figure}

Fig.\ref{fg5}(a) illustrates BPO/BPPs and DEPPs. 
As the number of DEPPs increases, the stability of the BPO becomes stronger. 
The stability is related to both the size of the basin of attraction \cite{ca3} \cite{ca4} and error-correction capability. 
Here we give several notations. 
Let $\bm{X}_- \equiv (0, \cdots, 0)$ and 
let $\bm{X}_+ \equiv (1, \cdots, 1)$. 
$\bm{X}_-$ is a fixed point and $\bm{X_+}$ is its DEPP: $F(\bm{X}_-)=\bm{X}_-$, $F(\bm{X}_+)=\bm{X}_-$. 
Except for these two end points, 
let $\bm{B}_e$ be a set of binary vectors containing an even number of 1s and 
let $\bm{B}_o$ be a set of binary vectors containing an odd number 1s. 
All the $N$-dimensional binary vectors are given by
$\bm{B}^N=\bm{B}_e\cup\bm{B}_o\cup \bm{X}_+ \cup \bm{X}_-$. 
Then we have
{\proposition 
Let $N \ge 5$ be an odd integer in PXCAs. 
The vector form function $F(B_e)=B_e$ is a bijection and $B_e$ consists of BPPs. 
$F(B_o)=B_e$ is a bijection and $B_o$ consists of DEPPs. 
The set $B_o$ consists of DEPPs, which are the inverse patterns of the BPPs. 
The numbers of BPPs and DEPPs are both $2^{N-1}-1$. 
}
{\remark 
The proof can be found in \cite{ncom25}. 
The proof uses the fact that, for odd $N \ge 5$, if a binary vector $\bm{X}$ contains an even number of 1s ($\bm{X} \in B_e$) then its inverse pattern contains an odd number of 1s ($\bar{\bm{X}} \in B_o$). 
However, it is not valid for even $N$. 
Since analysis in even $N$ is difficult, we focus on odd $N$ in this paper:}
\begin{equation}
\mbox{Dimension }N \in \{5, 7, 9, 11, 13, \cdots \}, \ 
\mbox{shift-part parameter } r \in \{1, \cdots, N \} 
\end{equation}
{\definition
A BPO is called a long binary periodic orbit (LBPO) if its period is $(2^N - 2)/2$: the LBPO consists of all elements of $\bm{B}_e$. 
All elements of $\bm{B}_o$ are DEPPs of the LBPO: stability is very strong. 
}  
Fig.~\ref{fg5} (b) shows an LBPO and its DEPPs for $N=7$ and $r=3$ .
{\remark
If a PXCA generates a LBPO, we can grasp all the binary vectors: $\bm{B}^N$ consists of BPPs of the LBPO, DEPPs of the LBPO, and two endpoints ($\bm{X}_-$ and $\bm{X}_+$). 
The LBPOs are different from M-sequences of the linear feedback shift register (LFSR \cite{lfsr1} \cite{lfsr2}) because the M-sequence has no DEPP and is not stable. 
If a PXCA does not have an LBPO, then it has multiple BPOs. 
Since it is difficult to clarify the number and period of all BPOs,   
we focus on the LBPO hereafter. 
}

\section{Theoretical Analysis}
We derive the transition matrix of PXCA over g$\mathrm{GF}(2)$. 
Using the characteristic polynomial of the transition matrix, we can analyze the LBPO algebraically. 
Since the XOR is addition over $\mathrm{GF}(2)$, the PXCA of Eq.~\eqref{pxca} can be represented by a dynamical system over $\mathrm{GF}(2)$. 
First,] Rule 90 is described by the ($N \times N$) transition matrix: 
\begin{equation}
\bm{R}= \bm{L}+\bm{L}^{-1}=
\begin{pmatrix}
0      & 1      & 0      & \cdots & 0      & 1      \\
1      & 0      & 1      & \ddots &        & 0      \\
0      & 1      & 0      & \ddots & \ddots & \vdots \\
\vdots & \ddots & \ddots & \ddots & 1      & 0      \\
0      &        & \ddots & 1      & 0      & 1      \\
1      & 0      & \cdots & 0      & 1      & 0
\end{pmatrix}
    \in \mathrm{GF}(2)^{(N \times N)}, \  \ 
\bm{X}^{t+1} = \bm{R} \bm{X}^t
\label{tmat90}
\end{equation}
where $\bm{L}$ is the cyclic shift matrix. 
The partial-shift permutation can be represented by the matrix
\begin{equation}
    \bm{P}_r=
    \begin{pmatrix}
        \bm{S}_r & \bm{0}_{r \times (N-r)} \\
        \bm{0}_{(N-r) \times r} & \bm{I}_{N-r}
    \end{pmatrix}
    \in \mathrm{GF}(2)^{(N \times N)}
    \label{pmat}
\end{equation}
where $\bm{I}_K$ is the $K$-dimensional identity and $\bm{S}_r$ is the $r \times r$ shift matrix: 
\begin{equation}
    \bm{S}_r=
    \begin{pmatrix}
        0 & 0 & 0 & \cdots & 1\\
        1 & 0 & 0 & \ddots & \vdots\\
        \vdots & \ddots & \ddots & \ddots & 0\\
        0 & \cdots & 1 & 0 & 0\\
        0 & 0 & \cdots & 1 & 0
    \end{pmatrix}\in \mathrm{GF}(2)^{(r\times r)}
    \end{equation}
Composing $\bm{P}_r$ and $\bm{R}$, we obtain the transition matrix of the PXCA:  
\begin{equation}
    \bm{A}=\bm{P}_r\bm{R}
    \in\mathrm{GF}(2)^{N\times N}, \ \ 
    \bm{X}^{t+1}=\bm{A}\bm{X}^t. 
    \label{tmat}
\end{equation}
The characteristic polynomial of the transition matrix $\bm{A}$ is 
\begin{equation}
    \chi_{\bm{A}}(x)=
    \det \left( x \bm{I}_N + \bm{A} \right)
    \in \mathrm{GF}(2)[x]
\end{equation}
Using this characteristic polynomial, we obtain Propositions 2 and 3.  
\medskip
\noindent
{\proposition 
The transition matrix $\bm{A}$ of Eq.~\eqref{tmat} has a zero eigenvalue thereby the characteristic polynomial has $x$ as a factor:  
\begin{equation}
    \chi_{\bm{A}}(x) = x p(x) \equiv q(x) \in \mathrm{GF}(2)[x]
    \label{primitive}
\end{equation}
}

\noindent
{\it Proof. }
Let $\bm{e}=(1, \cdots, 1 )^\top$ be the all-one vector. 
Applying the Rule 90 transition matrix $\bm{R}$ to $\bm{e}$, we notice
\begin{equation*}
\bm{R} \bm{e} = (\bm{L}+\bm{L}^{-1}) \bm{e} = \bm{L}\bm{e} + \bm{L}^{-1}\bm{e} = 
\bm{e} + \bm{e} = 
\bm{0}
\end{equation*}
Since $\bm{A}$ is the composition of the partial-shift permutation matrix $\bm{P}_r$ and $\bm{R}$, we obtain
\begin{equation*}
    \bm{A}\bm{e}
    =
    \bm{P}_r\bm{R}\bm{e}
    =
    \bm{P}_r\bm{0}
    =
    \bm{0}
\end{equation*}
That is, $\bm{A}$ has zero 
eigenvalue and $\bm{e}$ is its eigenvector. Hence
\begin{equation*}
\chi_{\bm{A}}(0) = 
\det ( \bm{A} ) = 0.
\end{equation*}
Thus, the characteristic polynomial has $x$ as a factor:  
\begin{equation*}
\chi_{\bm{A}}(x)= x p(x) 
  \hspace{10mm}  \square
\end{equation*}
{\definition
A characteristic polynomial $q(x) = x p(x)$, where $p(x)$ is primitive, is called {\it a half-primitive polynomial}.
$p(x)$ is called {\it a primitive factor.}}

{\proposition
If $q(x)$ is a half-primitive polynomial over $\mathrm{GF}(2)$, then the PXCA generates an LBPO with period $2^{N-1}-1$.
}

\noindent
{\it Proof.} 
Since $p(x)$ is a primitive factor, the half-primitive polynomial $q(x)=xp(x)$ guarantees the generation of a BPO with period $2^{N-1}-1$ \cite{lfsr1} \cite{lfsr2}. 
According to Proposition 1, all elements of $\bm{B}_e$ are BPPs, and their total number is $2^{N-1}-1$. 
Therefore, the elements of $\bm{B}_e$ constitute the BPO with period $2^{N-1}-1$. 
On the other hand, all elements of $\bm{B}_o$ are DEPPs and fall directly into the BPO. 
Consequently, the resulting BPO is an LBPO. 
$\ \ \square$
%


\noindent
Here we show two examples. 

\noindent
{\bf Example 1.} Fig.~\ref{fg3} (b) shows a PXCA for $N=9$ and $r=5$. 
For $N=9$, Rule 90 transition matrix in Eq.~\eqref{tmat90} is
\begin{equation*}
    \bm{R}=
    \begin{pmatrix}
        0&1&0&0&0&0&0&0&1\\
        1&0&1&0&0&0&0&0&0\\
        0&1&0&1&0&0&0&0&0\\
        0&0&1&0&1&0&0&0&0\\
        0&0&0&1&0&1&0&0&0\\
        0&0&0&0&1&0&1&0&0\\
        0&0&0&0&0&1&0&1&0\\
        0&0&0&0&0&0&1&0&1\\
        1&0&0&0&0&0&0&1&0
    \end{pmatrix}
\end{equation*}
For $r=5$, the partial-shift permutation matrix in Eq.~\eqref{pmat} is
\begin{equation*}
    \bm{P}_5=
    \begin{pmatrix}
        0&0&0&0&1&0&0&0&0\\
        1&0&0&0&0&0&0&0&0\\
        0&1&0&0&0&0&0&0&0\\
        0&0&1&0&0&0&0&0&0\\
        0&0&0&1&0&0&0&0&0\\
        0&0&0&0&0&1&0&0&0\\
        0&0&0&0&0&0&1&0&0\\
        0&0&0&0&0&0&0&1&0\\
        0&0&0&0&0&0&0&0&1
    \end{pmatrix}
\end{equation*}
%
%
The characteristic polynomial of the transition matrix $\bm{A}_5 = \bm{P}_5 \bm{R}$ is
\begin{equation*}
        \chi_{\bm{A}_5}(x) =
        \det\left(x\bm{I}_9-\bm{A}_5\right) =
        xp(x), \ p(x)=x^8+x^6+x^3+x^2+1
\end{equation*}
The $p(x)$ is a primitive factor of degree 8 over $\mathrm{GF}(2)$.  
Hence, the PXCA generates the LBPO with period $2^8-1$, as shown in Fig.~\ref{fg4} (b). 
All $2^8-1$ elements in $\bm{B}_o$ are DEPPs and fall directly into the LBPO. 

\noindent
{\bf Example 2.}  Fig.~\ref{fg3}(c) shows a PXCA for $N=9$ and $r=7$. 
The characteristic polynomial of the transition matrix 
$\bm{A}_7 = \bm{P}_7\bm{R}$ is
\begin{equation*}
        \chi_{\bm{A}_7}(x) =
        \det\left(x\bm{I}_9-\bm{A}_7\right) =
        xp(x), \ p(x)=x^8+x^5+x^3+x+1
\end{equation*}
The polynomial $p(x)$ is a primitive factor of degree 8 over $\mathrm{GF}(2)$ 
and the PXCA generates the LBPO with period $2^8-1$, as shown in Fig.~\ref{fg4}(c). 
All $2^8-1$ elements in $\bm{B}_0$ are DEPPs and fall directly into the LBPO. 
Based on the above discussions, the problem is  
\begin{itemize}
\item{\it 
exploring the values of the shift-part parameter $r$ and the odd-dimension $N$ that enable the generation of LBPOs.}
\end{itemize}
In order to solve this problem, we present a simple algorithm in Algorithm~\ref{alg1}. 
In the algorithm, characteristic polynomial is computed using SymPy and the primitivity of $p(x)$ is tested using the \texttt{galois} Python package. 
Applying the algorithm to PXCAs for $N \in \{5, \cdots, 203\}$, we give the main results summarized in Table~\ref{tb1} for 
$N \in \{5, \cdots, 51 \}$ and Table \ref{tb2} for $N \in \{53, \cdots, 203 \}$. 
Table~\ref{tb1} shows the values of $r$ and $N$ that enable the generation of LBPOs characterized by the primitive factor $p(x)$ in the table. 
Table~\ref{tb2} shows the values of $r$ and $N$ that enable the generation of LBPOs. 
For convenience in the table space, primitive factors $p(x)$ are omitted. 
We can see that $r=1$ and $r=N$ do not appear in the Tables: the Rule 90 and the full shift permutation do not enable the generation of LBPOs  (except for $r=N=5$). 
The PXCA does not generate LBPOs for
\[
N \in \{ 13, 25, 35, 37, 41, 45, 49, 65, 79, 113, 133, 143, 151, 153, 161, 169 \}
\]
{\remark 
Calculation for $N > 203$ is theoretically possible; however, the period $2^{203-1}-1$ is already sufficiently long for  practical purposes. 
As stated earlier, our previous study \cite{ncom25} investigated a PXCA with full-shift permutation ($r=N$) and suggested the generation of various BPOs. 
However, it addressed neither partial-shift connections, LBPOs, nor the characteristic-polynomial analysis presented in this section. 
}
\begin{algorithm}
\caption{Exploring $N$ and $r$ for primitive factor $p(x)$: LBPO.}
\label{alg1}
\begin{algorithmic}

\STATE $N = $ odd \hfill // dimension
\STATE $r \in \{1,\ldots,N\}$ \hfill // shift-part parameter

\STATE Construct the rule-90 matrix $\bm{R}$ \hfill

\FOR{$r = 1$ to $N$}
\STATE Construct the partial-shift permutation matrix $\bm{P}_r$

\STATE $\bm{A} \gets \bm{P}_r\bm{R}$
\hfill // transition matrix of PXCA

\STATE $\chi_{\bm{A}}(x)
\gets \det(x\bm{I}_N-\bm{A})$
\hfill // with \texttt{sympy.Matrix.charpoly()}

\STATE $p(x) \gets \chi_{\bm{A}}(x)/x$

\IF {$p(x)$ is primitive over $\mathrm{GF}(2)$ (with \texttt{galois.poly())}}

    \STATE  Record $N$, $r$, and $p(x)$ 
    \hfill // LBPO generation
\ENDIF

\ENDFOR
\end{algorithmic}
\end{algorithm}
\begin{table}[ ]
\centering
{
\footnotesize
\setlength{\tabcolsep}{4pt}
\renewcommand{\arraystretch}{1.2}
\caption{Dimension $N$, parameter $r$, and primitive factors of LBPOs. }
\label{tb1}
\begin{tabular}{c|c|>{\raggedright\arraybackslash}p{0.76\textwidth}}
\hline
$N$ & $r$ & Primitive factor $p(x)$ \\
\hline
5 & 3 & $x^{4}+\allowbreak x+\allowbreak 1$ \\
 & 5 & $x^{4}+\allowbreak x^{3}+\allowbreak 1$ \\
\hline
7 & 3 & $x^{6}+\allowbreak x+\allowbreak 1$ \\
\hline
9 & 5 & $x^{8}+\allowbreak x^{6}+\allowbreak x^{3}+\allowbreak x^{2}+\allowbreak 1$ \\
 & 7 & $x^{8}+\allowbreak x^{5}+\allowbreak x^{3}+\allowbreak x+\allowbreak 1$ \\
\hline
11 & 9 & $x^{10}+\allowbreak x^{8}+\allowbreak x^{7}+\allowbreak x^{2}+\allowbreak 1$ \\
\hline
15 & 3 & $x^{14}+\allowbreak x^{8}+\allowbreak x^{4}+\allowbreak x+\allowbreak 1$ \\
 & 7 & $x^{14}+\allowbreak x^{10}+\allowbreak x^{8}+\allowbreak x^{6}+\allowbreak x^{5}+\allowbreak x^{3}+\allowbreak x^{2}+\allowbreak x+\allowbreak 1$ \\
 & 8 & $x^{14}+\allowbreak x^{13}+\allowbreak x^{12}+\allowbreak x^{11}+\allowbreak x^{9}+\allowbreak x^{8}+\allowbreak x^{6}+\allowbreak x^{3}+\allowbreak x^{2}+\allowbreak x+\allowbreak 1$ \\
 & 9 & $x^{14}+\allowbreak x^{12}+\allowbreak x^{8}+\allowbreak x^{7}+\allowbreak x^{6}+\allowbreak x^{4}+\allowbreak 1$ \\
 & 12 & $x^{14}+\allowbreak x^{13}+\allowbreak x^{12}+\allowbreak x^{11}+\allowbreak x^{10}+\allowbreak x^{9}+\allowbreak x^{7}+\allowbreak x^{5}+\allowbreak x^{3}+\allowbreak x+\allowbreak 1$ \\
 & 14 & $x^{14}+\allowbreak x^{13}+\allowbreak x^{12}+\allowbreak x^{11}+\allowbreak x^{10}+\allowbreak x^{9}+\allowbreak x^{8}+\allowbreak x^{6}+\allowbreak x^{4}+\allowbreak x^{2}+\allowbreak 1$ \\
\hline
17 & 3 & $x^{16}+\allowbreak x^{12}+\allowbreak x^{10}+\allowbreak x^{8}+\allowbreak x^{2}+\allowbreak x+\allowbreak 1$ \\
 & 8 & $x^{16}+\allowbreak x^{15}+\allowbreak x^{14}+\allowbreak x^{13}+\allowbreak x^{12}+\allowbreak x^{9}+\allowbreak x^{8}+\allowbreak x^{7}+\allowbreak x^{5}+\allowbreak x^{3}+\allowbreak x^{2}+\allowbreak x+\allowbreak 1$ \\
\hline
19 & 3 & $x^{18}+\allowbreak x^{12}+\allowbreak x^{10}+\allowbreak x^{8}+\allowbreak x^{2}+\allowbreak x+\allowbreak 1$ \\
\hline
21 & 12 & $x^{20}+\allowbreak x^{19}+\allowbreak x^{18}+\allowbreak x^{17}+\allowbreak x^{16}+\allowbreak x^{15}+\allowbreak x^{13}+\allowbreak x^{12}+\allowbreak x^{11}+\allowbreak x^{9}+\allowbreak x^{7}+\allowbreak x^{3}+\allowbreak x^{2}+\allowbreak x+\allowbreak 1$ \\
\hline
23 & 9 & $x^{22}+\allowbreak x^{20}+\allowbreak x^{16}+\allowbreak x^{14}+\allowbreak x^{12}+\allowbreak x^{7}+\allowbreak 1$ \\
\hline
27 & 8 & $x^{26}+\allowbreak x^{25}+\allowbreak x^{24}+\allowbreak x^{23}+\allowbreak x^{21}+\allowbreak x^{20}+\allowbreak x^{18}+\allowbreak x^{14}+\allowbreak x^{10}+\allowbreak x^{7}+\allowbreak x^{6}+\allowbreak x+\allowbreak 1$ \\
 & 20 & $x^{26}+\allowbreak x^{25}+\allowbreak x^{24}+\allowbreak x^{23}+\allowbreak x^{18}+\allowbreak x^{17}+\allowbreak x^{15}+\allowbreak x^{14}+\allowbreak x^{13}+\allowbreak x^{12}+\allowbreak x^{10}+\allowbreak x^{9}+\allowbreak x^{7}+\allowbreak x^{6}+\allowbreak x^{5}+\allowbreak x^{4}+\allowbreak x^{2}+\allowbreak x+\allowbreak 1$ \\
\hline
29 & 14 & $x^{28}+\allowbreak x^{27}+\allowbreak x^{24}+\allowbreak x^{23}+\allowbreak x^{21}+\allowbreak x^{17}+\allowbreak x^{16}+\allowbreak x^{14}+\allowbreak x^{12}+\allowbreak x^{11}+\allowbreak x^{10}+\allowbreak x^{7}+\allowbreak x^{4}+\allowbreak x^{3}+\allowbreak 1$ \\
 & 15 & $x^{28}+\allowbreak x^{22}+\allowbreak x^{20}+\allowbreak x^{18}+\allowbreak x^{16}+\allowbreak x^{13}+\allowbreak x^{12}+\allowbreak x^{11}+\allowbreak x^{9}+\allowbreak x^{7}+\allowbreak x^{6}+\allowbreak x^{5}+\allowbreak x^{4}+\allowbreak x^{3}+\allowbreak x^{2}+\allowbreak x+\allowbreak 1$ \\
\hline
31 & 17 & $x^{30}+\allowbreak x^{28}+\allowbreak x^{24}+\allowbreak x^{16}+\allowbreak x^{15}+\allowbreak x^{14}+\allowbreak x^{12}+\allowbreak x^{8}+\allowbreak 1$ \\
 & 24 & $x^{30}+\allowbreak x^{29}+\allowbreak x^{28}+\allowbreak x^{27}+\allowbreak x^{26}+\allowbreak x^{25}+\allowbreak x^{24}+\allowbreak x^{23}+\allowbreak x^{22}+\allowbreak x^{21}+\allowbreak x^{20}+\allowbreak x^{19}+\allowbreak x^{17}+\allowbreak x^{16}+\allowbreak x^{14}+\allowbreak x^{11}+\allowbreak x^{9}+\allowbreak x^{8}+\allowbreak x^{6}+\allowbreak x^{3}+\allowbreak x^{2}+\allowbreak x+\allowbreak 1$ \\
\hline
33 & 14 & $x^{32}+\allowbreak x^{31}+\allowbreak x^{28}+\allowbreak x^{27}+\allowbreak x^{25}+\allowbreak x^{24}+\allowbreak x^{23}+\allowbreak x^{21}+\allowbreak x^{19}+\allowbreak x^{18}+\allowbreak x^{17}+\allowbreak x^{16}+\allowbreak x^{15}+\allowbreak x^{14}+\allowbreak x^{13}+\allowbreak x^{9}+\allowbreak x^{7}+\allowbreak x^{6}+\allowbreak x^{5}+\allowbreak x^{4}+\allowbreak x^{3}+\allowbreak x^{2}+\allowbreak 1$ \\
 & 30 & $x^{32}+\allowbreak x^{31}+\allowbreak x^{17}+\allowbreak x^{16}+\allowbreak x^{15}+\allowbreak x^{14}+\allowbreak x^{13}+\allowbreak x^{12}+\allowbreak x^{11}+\allowbreak x^{10}+\allowbreak x^{9}+\allowbreak x^{8}+\allowbreak x^{7}+\allowbreak x^{6}+\allowbreak x^{5}+\allowbreak x^{4}+\allowbreak x^{3}+\allowbreak x^{2}+\allowbreak 1$ \\
\hline
39 & 7 & $x^{38}+\allowbreak x^{34}+\allowbreak x^{32}+\allowbreak x^{30}+\allowbreak x^{26}+\allowbreak x^{24}+\allowbreak x^{22}+\allowbreak x^{20}+\allowbreak x^{18}+\allowbreak x^{16}+\allowbreak x^{6}+\allowbreak x^{5}+\allowbreak x^{4}+\allowbreak x^{3}+\allowbreak x^{2}+\allowbreak x+\allowbreak 1$ \\
 & 13 & $x^{38}+\allowbreak x^{36}+\allowbreak x^{34}+\allowbreak x^{28}+\allowbreak x^{24}+\allowbreak x^{22}+\allowbreak x^{18}+\allowbreak x^{14}+\allowbreak x^{11}+\allowbreak x^{10}+\allowbreak x^{7}+\allowbreak x^{6}+\allowbreak x^{3}+\allowbreak x^{2}+\allowbreak 1$ \\
 & 31 & $x^{38}+\allowbreak x^{34}+\allowbreak x^{32}+\allowbreak x^{29}+\allowbreak x^{27}+\allowbreak x^{25}+\allowbreak x^{23}+\allowbreak x^{21}+\allowbreak x^{19}+\allowbreak x^{17}+\allowbreak x^{15}+\allowbreak x^{13}+\allowbreak x^{11}+\allowbreak x^{9}+\allowbreak x^{7}+\allowbreak x^{6}+\allowbreak x^{5}+\allowbreak x^{3}+\allowbreak x^{2}+\allowbreak x+\allowbreak 1$ \\
\hline
43 & 9 & $x^{42}+\allowbreak x^{40}+\allowbreak x^{36}+\allowbreak x^{28}+\allowbreak x^{18}+\allowbreak x^{16}+\allowbreak x^{10}+\allowbreak x^{8}+\allowbreak x^{7}+\allowbreak x^{2}+\allowbreak 1$ \\
 & 13 & $x^{42}+\allowbreak x^{40}+\allowbreak x^{38}+\allowbreak x^{34}+\allowbreak x^{30}+\allowbreak x^{28}+\allowbreak x^{20}+\allowbreak x^{16}+\allowbreak x^{12}+\allowbreak x^{11}+\allowbreak x^{8}+\allowbreak x^{7}+\allowbreak x^{4}+\allowbreak x^{3}+\allowbreak 1$ \\
 & 26 & $x^{42}+\allowbreak x^{41}+\allowbreak x^{38}+\allowbreak x^{37}+\allowbreak x^{30}+\allowbreak x^{28}+\allowbreak x^{25}+\allowbreak x^{24}+\allowbreak x^{23}+\allowbreak x^{20}+\allowbreak x^{15}+\allowbreak x^{7}+\allowbreak x^{5}+\allowbreak x^{2}+\allowbreak 1$ \\
 & 35 & $x^{42}+\allowbreak x^{36}+\allowbreak x^{34}+\allowbreak x^{33}+\allowbreak x^{32}+\allowbreak x^{31}+\allowbreak x^{10}+\allowbreak x^{4}+\allowbreak x^{2}+\allowbreak x+\allowbreak 1$ \\
\hline
47 & 19 & $x^{46}+\allowbreak x^{40}+\allowbreak x^{36}+\allowbreak x^{32}+\allowbreak x^{30}+\allowbreak x^{24}+\allowbreak x^{17}+\allowbreak x^{15}+\allowbreak x^{4}+\allowbreak x+\allowbreak 1$ \\
\hline
51 & 14 & $x^{50}+\allowbreak x^{49}+\allowbreak x^{43}+\allowbreak x^{42}+\allowbreak x^{41}+\allowbreak x^{37}+\allowbreak x^{36}+\allowbreak x^{35}+\allowbreak x^{34}+\allowbreak x^{31}+\allowbreak x^{30}+\allowbreak x^{25}+\allowbreak x^{24}+\allowbreak x^{22}+\allowbreak x^{20}+\allowbreak x^{19}+\allowbreak x^{17}+\allowbreak x^{16}+\allowbreak x^{12}+\allowbreak x^{10}+\allowbreak x^{9}+\allowbreak x^{8}+\allowbreak x^{7}+\allowbreak x^{4}+\allowbreak 1$ \\
 & 16 & $x^{50}+\allowbreak x^{49}+\allowbreak x^{48}+\allowbreak x^{47}+\allowbreak x^{46}+\allowbreak x^{45}+\allowbreak x^{44}+\allowbreak x^{43}+\allowbreak x^{42}+\allowbreak x^{39}+\allowbreak x^{38}+\allowbreak x^{37}+\allowbreak x^{36}+\allowbreak x^{27}+\allowbreak x^{18}+\allowbreak x^{15}+\allowbreak x^{10}+\allowbreak x+\allowbreak 1$ \\
\hline
\end{tabular}
}
\end{table}
\clearpage

\begin{table}[]
\centering
\scriptsize
\setlength{\tabcolsep}{3pt}
\renewcommand{\arraystretch}{1}
\caption{Dimension $N$ and parameter $r$ for the generation of LBPOs}
\label{tb2}
\begin{tabular}{|c|>{\raggedright\arraybackslash}p{0.235\textwidth}|c|>{\raggedright\arraybackslash}p{0.235\textwidth}|c|>{\raggedright\arraybackslash}p{0.235\textwidth}|}
\hline
$N$ & $r$ & $N$ & $r$ & $N$ & $r$ \\
\hline
53 & 22, 46 & 103 & 55, 56 & 157 & 37, 59 \\
\hline
55 & 38 & 105 & 38, 101 & 159 & 13, 45, 125, 128 \\
\hline
57 & 25, 31, 45, 49 & 107 & 8, 15, 16, 17, 49, 50, 86 & 163 & 95, 144 \\
\hline
59 & 8, 19, 31, 37 & 109 & 11, 35, 61 & 165 & 14, 79, 80, 131 \\
\hline
61 & 27 & 111 & 12, 75 & 167 & 47, 83, 158 \\
\hline
63 & 9 & 115 & 74 & 171 & 150 \\
\hline
67 & 17, 31, 62 & 117 & 26, 83 & 173 & 27, 102, 136 \\
\hline
69 & 48, 59 & 119 & 37, 60, 65 & 175 & 79, 80 \\
\hline
71 & 23, 39, 43 & 121 & 53 & 177 & 80 \\
\hline
75 & 3, 8, 17, 28, 37, 38 & 123 & 13, 28, 46, 59, 87, 90 & 179 & 151, 164 \\
\hline
77 & 24 & 125 & 63, 68 & 181 & 124 \\
\hline
81 & 19, 53, 71, 78 & 127 & 56, 93, 98, 120 & 183 & 16, 26, 51, 111, 132 \\
\hline
83 & 42, 54 & 129 & 14, 28 & 185 & 134 \\
\hline
85 & 5, 8, 14, 80 & 131 & 129 & 187 & 29, 33, 109, 167 \\
\hline
87 & 31, 80, 84 & 135 & 42 & 189 & 20, 60 \\
\hline
89 & 32, 52, 59 & 137 & 8, 106, 108 & 191 & 91, 105 \\
\hline
91 & 73 & 139 & 98, 122 & 193 & 55 \\
\hline
93 & 14, 37, 76, 91 & 141 & 124 & 195 & 3 \\
\hline
95 & 33, 41 & 145 & 51 & 197 & 51, 55, 68 \\
\hline
97 & 19, 92 & 147 & 138 & 199 & 59 \\
\hline
99 & 28, 32, 75, 88 & 149 & 58, 74 & 201 & 134, 169 \\
\hline
101 & 19, 48, 91 & 155 & 69, 138 & 203 & 9, 103 \\
\hline
\end{tabular}
\end{table}

\section{FPGA-based hardware}
We present a simple FPGA-based hardware implementation of the PXCAs. 
The FPGA is reconfigurable digital hardware enabling high-speed processing over software based implementation in laptops \cite{rc4}. 
Both the partial-shift connection and the XOR rule facilitate efficient hardware implementation.
The design outline is shown in Algorithm \ref{alg2}, the SystemVerilog pseudocode.  
Giving the dimension $N$ and the shift-part parameter $r$, the PXCA is implemented on an FPGA board.  
\begin{algorithm}
    \caption{FPGA-based hardware design of PXCA}
    \label{alg2}
    \begin{algorithmic}
        \STATE $N = $odd \hfill // dimension
        \STATE $r \in \{1, \cdots, N \}, P[1:N]$ \hfill //  parameters
        \STATE \textbf{input} $ i_0[1:N]$
        \STATE \textbf{output} $x^{t}[1:N]$
        \STATE \textbf{reg} $y^{t}[1:N]$
        

        
        \FOR{$i = 1$ to $N$}
                \IF{$i==1$}
                    \STATE $P[i] \gets r$
                \ELSIF{$2 \le i \le r$}
                    \STATE $P[i] \gets i-1$
                \ELSE
                    \STATE $P[i] \gets i$
                \ENDIF
        \ENDFOR
            
        \FOR{$j = 1$ to $N$}
            \STATE $y^{t}[j] \gets x^{t}[(j-1)] \oplus x^{t}[(j+1)]$ \hfill // XOR
        \ENDFOR

        \FOR{$k = 1$ to $N$}
            \STATE $x^t[k] \gets y^{t}[P[k]]$ \hfill //permutation
        \ENDFOR

    \end{algorithmic}
\end{algorithm}

\begin{figure}[htb]
\centering
\includegraphics[width=0.6\columnwidth]{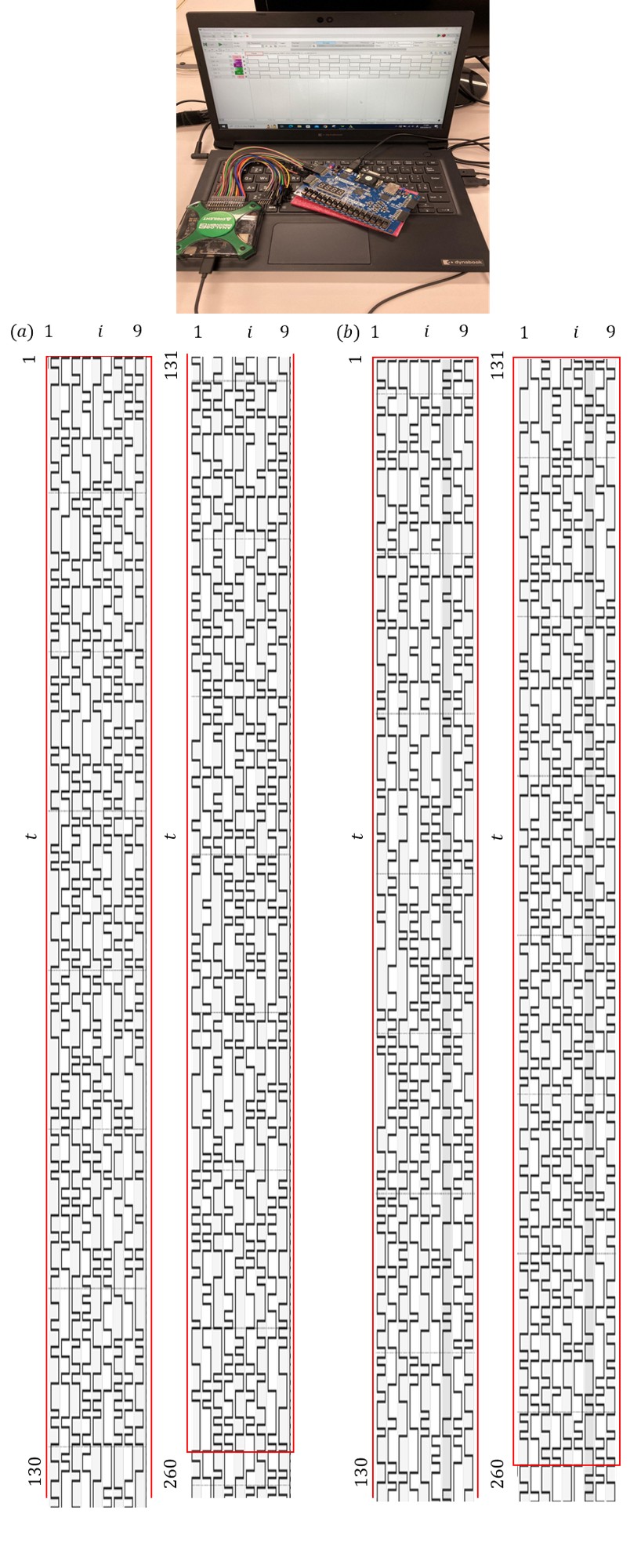}
\caption{Measured LBPOs for $N=9$ (period $255$) on the FPGA board.
(a) $r=5$ (as is Fig. \ref{fg4} (b)).
(b) $r=7$ (as is Fig. \ref{fg4} (c)).
}
\label{fg6}
\end{figure}

\noindent
The actual design/measurement environment is
\begin{itemize}
    \item Design: Vivado 2023.1 platform (Xilinx)
    \item FPGA: Xilinx Artix-7 XC7A35T-1CPG236C 
    \item Clock: 100 [kHz] (the default frequency 100 [MHz] is divided)
    \item Measurement: ANALOG DISCOVERY2, Waveforms 2015
\end{itemize}
Note that the hardware performance (e.g., speed and scalability) depends on the FPGA specification of rapid progress. 
The FPGA-based hardware converts LBPOs into electrical signals as shown in Fig~\ref{fg6}. 
The signals can be used in engineering application systems such as pseudo random number generators \cite{rnd} and switching circuits controllers \cite{pe1}. 
Owing to the strong stability of LBPOs, 
the electrical signals exhibit error-correction capability and robustness against external perturbations.

\section{Conclusions}
LBPOs generated by PXCAs are investigated in this paper. 
The PXCA is characterized by a partial-shift permutation and Rule 90. 
The PXCA can generate a wide variety of BPOs; in particular, we focus on LBPOs with strong stability.
By deriving characteristic polynomials over $\mathrm{GF}(2)$ from the transition matrix of the PXCA, we identify the values of the shift-part parameter $r$ and the odd dimension $N$ that enable the generation of LBPOs. 
Furthermore, using a simple FPGA-based hardware implementation, we obtain electrical signals corresponding to LBPOs for potential engineering applications. 

Many issues remain for future study, including the following:
1) Analysis of the generation mechanism of LBPOs.
2) Investigation of the relationship between PXCAs and LFSRs, which has not yet been fully clarified. 
In particular, the characteristic polynomial may provide a bridge between the PXCA and LFSR frameworks.
3) Exploration of engineering applications of various LBPOs.

\clearpage

\clearpage

\medskip


\begin{thebibliography}{99}

\bibitem{ca1}
S. Wolfram,
Cellular Automata and Complexity: Collected Papers, CRC Press, 2018.

\bibitem{ca2}
S. Wolfram, 
Computation theory of cellular automata, 
Commun. Math. Phys. 96 (1984) 15-57 

\bibitem{ca3}
A. Wuensche, M. Lesser, Global dynamics of cellular automata: an atlas of basin of attraction fields of one-dimensional cellular automata (Vol. 1),  Santa Fe Institute Press, 1992.

\bibitem{ca4}
J. E. Hanson and J. P. Crutchfield, 
The attractor-basin portrait of a cellular automaton,  
Journal of statistical physics, 66 (1992) 1415-1462 

\bibitem{ca5}
L. O. Chua, A nonlinear dynamics perspective of Wolfram's new kind of science, I, World Scientific, 2006.

\bibitem{ca6}
M. Sch\"{u}le, R. Stoop, 
A full computation-relevant topological dynamics classification of elementary cellular automata. 
Chaos, 22 (2012) 043143

\bibitem{rc1}
O. Yilmaz,  
Symbolic computation using cellular automata-based hyperdimensional computing, 
Neural Computation, 27 (2015) 2661-2692. 

\bibitem{rc2}
S. Nichele and M. S. Gundersen, Reservoir Computing Using
Nonuniform Binary Cellular Automata, Complex Systems, 26, 3  (2017) 225–245. 


\bibitem{rc3}
N. Babson,  C. Teusche, 
Reservoir computing with complex cellular automata. 
Complex Systems, 28 (2019) 433-455. 

\bibitem{rc4}
G. Tanaka, T. Yamane, J. B. H\'eroux, R. Nakane, N. Kanazawa, S. Takeda, H. Numata, D. Nakano, A. Hirose, 
Recent advances in physical reservoir computing: a review,  
Neural Networks, 115 (2019) 100-123.

\bibitem{error}
D. R. Chowdhury, S. Basu, I. S. Gupta, P. P. Chaudhuri, 
Design of CAECC - cellular automata based error correcting code, 
IEEE Trans. Comput. 43 (6): (1994) 759-764. 

\bibitem{data}
W. Wada, J. Kuroiwa, S. Nara, 
Completely reproducible description of digital sound data with cellular automata, 
Phys. Lett. A 306 (2002) 110-115. 

\bibitem{iconip22}
T. Okano and T. Saito, 
Permutation elementary cellular automata: analysis and application of simple examples, 
M. Tanveer et al. (Eds.): ICONIP 2022, LNCS 13623 (2023) 321-330

\bibitem{iconip24} 
M. Onuki and T. Saito, Analysis of direct stable binary periodic orbits in permutation elementary cellular automata, 
M. Mahmud et al. (Eds.): ICONIP 2024, CCIS 2282 (2025) 380-392

\bibitem{ncom25}
M. Onuki, Y. Suzuki and T. Saito,  
Permutation XOR cellular automata and direct stable periodic orbits, 
Neurocomputing 656 (2025) 131510

\bibitem{taka}
H. Udagawa, T. Okano and T. Saito, 
Permutation binary neural networks: analysis of periodic orbits and its applications, 
Discrete Contin. Dyn. Syst. Ser. B, 28, 1, (2023) 748-762. 

\bibitem{mikito}
M. Onuki, K. Saka and T. Saito, 
A variety of globally stable periodic orbits in permutation binary neural networks, 
Discrete Contin. Dyn. Syst. Ser. B, doi/10.3934/dcdsb.2023078  (2023)

\bibitem{lfsr1}
S. W. Golomb, Shift register sequences, revised 2nd edition, Aegean Park Press, 1982. 

\bibitem{lfsr2}
Z. Chang, M. F. Ezerman, S. Ling, H. Wang, 
The cycle structure of LFSR with arbitrary characteristic polynomial over finite fields, 
Cryptography and Communications, 10 (2018) 1183-1202. 



\bibitem{rnd}
A. Poorghanad, A. Sadr, A. Kashanipour, 
Generating high quality pseudo random number using evolutionary methods, 
in: Proc. IEEE/CIS, 9 (2008) 331-335. 

\bibitem{pe1}
W. Holderbaum, 
Application of neural network to hybrid systems with binary inputs.  
IEEE Trans. Neural Netw. 18, 4 (2007) 1254-1261

\bibitem{pe2}
P. W. Wheeler, J. Rodriguez, J. C. Clare, L. Empringham and A. Weinstein, 
Matrix converters: a technology review,  
IEEE Tran. Ind. Electron., 49, 2 (2002) 276-288

\end{thebibliography}
\end{document}